\input epsf
\magnification=\magstep1

%\raggedbottom

\hsize=6.35truein \vsize=9truein

\font\title=cmssdc10 at 15pt

\def\la{\lambda}

\def\ZZ{\bf Z}
\def\frame #1{\vbox{\hrule height.1pt
\hbox{\vrule width.1pt\kern 10pt
\vbox{\kern 10pt
\vbox{\hsize 3.5in\noindent #1}
\kern 10pt}
\kern 10pt\vrule width.1pt}
\hrule height0pt depth.1pt}}
\def\bframe #1{\vbox{\hrule height.1pt
\hbox{\vrule width.1pt\kern 10pt
\vbox{\kern 10pt
\vbox{\hsize 5in\noindent #1}
\kern 10pt}
\kern 10pt\vrule width.1pt}
\hrule height0pt depth.1pt}}
\baselineskip 15pt

\centerline{\title Will the swine strain crowd out the seasonal influenza strain?}
\bigskip

Rinaldo B. Schinazi

University of Colorado at Colorado Springs

Colorado Springs CO80933-7150 USA

rschinaz@uccs.edu

\footnote{}{Key words and phrases: competition models, stochastic process, influenza, swine strain}

\bigskip
\bigskip
{\bf Abstract.}  We use spatial and non spatial models to argue that competition alone may explain why  two influenza strains do not usually coexist. The more virulent strain is likely to crowd out the less virulent one. This can be seen as a consequence of the Exclusion Principle of Ecology. We exhibit, however, a spatial model for which coexistence is possible.
\bigskip
{\bf 1. Introduction.}  The seasonal flu strain has been a lot less prevalent during the 2009/2010 influenza season than during the previous years, see Fluview (the weekly CDC inluenza report) . On the other hand, some time during Spring 2009 the new so called swine strain has appeared. There  seems to be a relation between these two events and at this point in time (March 2010) it seems likely that the swine strain will crowd out the seasonal strain. In this paper we propose to explain this phenomenon using competition models.  We will use spatial and non spatial models to show that coexistence of two strains is unlikely due to competition alone. We will show along the way that the topology of space may be crucial.

A competing  explanation of the non coexistence of the two influenza strains is cross immunity. For instance, immunity may explain why older generations have not been as much affected as the younger ones in this swine epidemic. It may be due to some previous exposure to a similar strain, see the Discussion in Greenbaum et al. (2009). However, using a cross immunity argument to explain why the swine strain crowds out the seasonal one may be more difficult. The hypothesis would be that the swine strain must confer some  immunity against the seasonal flu. But, clearly the seasonal strain does not confer any immunity to the swine strain: after all even young people (the group most severely affected by the swine strain) have usually been exposed to the seasonal strain and do not seem to be protected against the swine strain. Hence, for this argument to work the swine strain must confer some immunity against the seasonal strain but the seasonal strain cannot confer any immunity against the swine strain. In contrast to this cross immunity hypothesis we argue in this paper that even in models for which there is no immunity at all (every individual that recovers is immediately susceptible again!) coexistence of two competing strains is rather unlikely.
\bigskip
{\bf 2. The ODE Model.}  Our first model is a system of ordinary differential equations. Let $u_1(t)$ and $u_2(t)$  be the density of individuals infected at time $t$ with strains 1 and 2, respectively .  We set
$$u_1'=\la_1u_1u_0-\delta_1u_1$$
$$u_2'=\la_2u_2u_0-\delta_2u_2$$
where $u_0(t)$ is the density of susceptible individuals at time $t$. In words,  individuals infected with strain $i$ infect susceptible individuals at rate $\la_i$ and get healthy at rate $\delta_i$, for $i=1,2$.  Since
$u_0(t)+u_1(t)+u_2(t)=1$ as soon as an infected individual gets healthy it is back in the susceptible pool. 

Let 1 be the seasonal and 2 be the swine strains. Early reports indicate that the swine strain may be more virulent than the seasonal strain, see Fraser et al. (2009). Under that assumption, 
$${\la_1\over\delta_1}<{\la_2\over\delta_2}.$$ 

Assume also that at some point in time the ODE model is at the equilibrium $(0, 1-{\delta_2\over\la_2})$. That is, there is no seasonal strain and the swine strain is in equilibrium. Now introduce a little bit of seasonal strain (small $u_1$). Will the seasonal strain be able to grow? Using that $u_1$ is almost 0 and that $u_2$ is almost $1-{\delta_2\over\la_2}$ we make the approximation
$$u_0=1-u_1-u_2\sim 1-(1-{\delta_2\over\la_2})={\delta_2\over\la_2}.$$
Hence,
$$u_1'\sim \la_1u_1{\delta_2\over\la_2}-\delta_1 u_1=u_1(\la_1{\delta_2\over\la_2}-\delta_1).$$
Since we are assuming that ${\la_1\over\delta_1}<{\la_2\over\delta_2}$ we get $u'_1<0$. That is, under these assumptions and according to this model the seasonal flu will not take hold.

In fact this system of ODE is a particular case of a well-known competition model. For the general version of this model it is known that one of the strains will vanish, see Exercise 3.3.5 in Hofbauer and Sigmund (1998). The point is that we have two populations (the population of individuals infected with strain 1 and the population of individuals infected with strain 2) that compete for a single resource (the susceptible individuals). It turns out that in such a model one population will drive the other one out. This is a particular case of the so called "Exclusion Principle" of Ecology: if the number of populations is larger than the number of resources all the populations cannot subsist in the long run, see 5.4 in Hofbauer and Sigmund (1998).
\bigskip
{\bf 2. The spatial stochastic model.} In the preceding model there is no space structure and all individuals can somehow be seen as neighbors. In this section we go the other extreme where there is a rigid space structure and each individual has a fixed number of neighbors. 

We now describe the multitype contact process, see Neu\-hau\-ser (1992).
Let $S$ be the integer lattice $\ZZ^d$ ( $d$ is the dimension) or the homogeneous tree ${\bf T}_d$ for which each site has $d+1$ neighbors.  The
system is described by a configuration $\xi\in \{0,1,2\}^S$,
where $\xi(x) =0$ means that site $x$ is occupied by a susceptible individual, 
$\xi(x)=1$ means that $x$ is occupied by an individual infected by strain
$1$ and $\xi(x)=2$ means that $x$ is occupied by an individual infected by strain
$2$.   If $S$ is $\ZZ^d$ then each site has $2d$ neighbors, if $S$ is ${\bf T}_d$ then each site has $d+1$ neighbors.
For $x\in S$ and $\xi\in\{0,1,2\}^{S}$, let
$n_1(x,\xi)$ and $n_2(x,\xi)$ denote the number of neighbors of $x$ that are infected by strain 1 and strain 2, respectively.

The multitype contact process $\xi_t$ with birth rates
$\lambda_1,\lambda_2$ makes  transitions at $x$ in configuration $\xi$
$$\eqalign{
 1\to 0 &\hbox{ at rate 1 }\cr
 2\to 0 &\hbox{ at rate 1 }\cr
0\to 1 &\hbox{ at rate }\lambda_1 n_1(x,\xi),\cr
0\to 2 &\hbox{ at rate }\lambda_2 n_2(x,\xi),\cr}
$$
In words, a susceptible individual gets infected by an infected neighbor at rates $\la_1$ or $\la_2$, depending on which strain the neighbor is infected with. An infected individual gets healthy (and is immediately susceptible again) at rate 1. Note that compared to the ODE model we are assuming in this model that $\delta_1=\delta_2=1$. This is so because most of the mathematical results have been proved under the assumption $\delta_1=\delta_2$. We take this common value to be 1 to minimize the number of parameters.

The multitype contact process is a generalization of the basic contact process which has only one type. 
The transitions of the basic contact process are given by
 $$\eqalign{
 1\to 0 &\hbox{ at rate 1 }\cr
0\to 1 &\hbox{ at rate }\lambda_1 n_1(x,\xi),\cr}
$$
There exists a critical value $\lambda_c$ (whose exact value is not known and which depends on the graph the model is on) for this model. If $\lambda_1>\lambda_c$ then starting with even a single infected individual there is a positive probability of having infected individuals at all times somewhere in the graph. On the other hand if $\lambda_1\leq \lambda_c$ then starting from any finite number of infected individuals all the infected individuals will disappear after a finite time. See Liggett (1999) for more on the basic contact process on the square lattice and on trees.

\medskip
{\bf 2.1 The space is the square lattice $\ZZ^d$.}  
We now go back to the multitype contact process. Assume that $\la_2>\la_c$ and $\la_2>\la_1$ then there is no coexistence of strains 1 and 2 in the sense that
$$\lim_{t\to\infty} P(\xi_t(x)=1,\xi_t(y)=2)=0$$
for any sites $x$ and $y$ in $\ZZ^d$. In fact, strain 2 always drives out strain 1 in the following sense. Conditioned on strain 2 not disappearing then
$$\lim_{t\to\infty} P(\xi_t(x)=1)=0,$$
for any site $x$ in $\ZZ^d$
and any initial configuration. See Theorem 2 in Cox and Schinazi (2009) and also Neuhauser (1992). 
Hence, assuming that $\la_2>\la_1$ (that is, strain 2 is more virulent than strain 1) this model too predicts that the seasonal flu will be crowded out by the swine strain. The spatial structure seems to have no influence on the outcome. The next section will show that this is not always so and and that a different (more crowded) space structure allows coexistence.
\medskip
{\bf 2.2 The space is the tree ${\bf T_d}$.}  There is a fundamental difference between the basic contact process on the square lattice and the same model on the tree. There are two (instead of one) critical values for the basic contact process on the tree. The definition of $\la_c$ is as  before. We also define another critical value $\la_{cc}$ in the following way. Consider the basic (one type) contact process with birth rate $\la_1$. Let $O$ be a fixed site on the tree or square lattice.  There is a positive probability that site $O$ will be infected infinitely often if and only if $\la_1>\la_{cc}$, starting with at least one infected individual. It turns out that $\la_c<\la_{cc}$ on the tree but $\la_c=\la_{cc}$ on the square lattice.  

We have the following result for the multitype contact process on the tree. Let $\la_1$ and $\la_2$ be in $(\la_c,\la_{cc})$ then strains 1 and 2 may coexist on the tree in the sense that there is a positive probability of having individuals infected with strain 1 and others with strain 2 at all times. See Theorem 1 in Cox and Schinazi (2009). Note that coexistence occurs even for $\la_1<\la_2$ but both parameters need to be in the rather narrow interval $(\la_c,\la_{cc})$. This result shows that space structure may be crucial in allowing coexistence.
\bigskip
{\bf 3. Discussion.} Our models show that at least in theory coexistence of two competing strains is unlikely. Coexistence is however possible for the multitype contact process on a tree. The tree can be thought of as a model for high density populations (in a ball of radius $r$ there are $(d+1)d^{r-1}$ individuals on  the tree ${\bf T_d}$ but only about $r^d$ on the lattice ${\ZZ^d}$). In order to have coexistence both infection rates cannot be too low or too high but may be unequal. In all other cases there will be no coexistence on the tree and there is never coexistence on $\ZZ^d$ unless $\la_1$ is exactly equal to $\la_2$, a rather unlikely possibility, see Neuhauser (1992). Interestingly the behavior of the ODE model (or mean-field) is the same as the behavior of the model on $\ZZ^d$ but different from the model on the tree. In general, it is expected that the model on the tree be closer to the mean-field model than to the model on $\ZZ^d$. This is not so in this example.

In all our models all infected individuals are put back in the susceptible pool as soon as they recover.   This is admittedly not realistic. However, it shows that there are other possible explanations of non coexistence than cross immunity arguments.
\bigbreak
{\bf References.}

\medskip\noindent
J.T. Cox and R.B.Schinazi (2009). Survival and coexistence for a multitype contact process. Annals of probability {\bf 37}, 853-876

\medskip\noindent
C. Fraser et al. (2009).Pandemic Potential of a Strain of Influenza A (H1N1): Early Findings. Science,  {\bf 324}, 1557 - 1561

\medskip\noindent
J.A. Greenbaum et al. (2009). Pre-existing immunity against swine-origin H1N1 influenza viruses in the general human population. PNAS, {\bf 106}, 20365-20370.

\medskip\noindent
T.Liggett (1999). {\it Stochastic interacting systems: contact, voter and
exclusion processes,} Springer-Verlag, Berlin.

\medskip\noindent 
C. Neuhauser (1992). Ergodic theorems for the multitype contact
process.  Probab. Theory Related Fields, {\bf 91}, 467-506.

\bye